\documentclass{article}

\usepackage{amssymb}
\usepackage{latexsym}
\usepackage{extarrows}
\usepackage[dvips]{graphicx}

\title{ Dense  $3$-uniform hypergraphs containing a large clique}
\author{Biao Wu \thanks{ College of Mathematics and Econometrics, Hunan University, Changsha 410082, P.R. China. Email: wubiao@hnu.edu.cn.} \and Yuejian Peng \thanks{Corresponding author. Institute of Mathematics, Hunan University, Changsha 410082, P.R. China. Email: ypeng1@hnu.edu.cn. Partially  supported by National Natural Science Foundation of China (No. 11271116).} }
\date{}

\newtheorem{theo}{Theorem}[section]

\newtheorem{remark}[theo]{Remark}
\newtheorem{defi}[theo]{Definition}

\newtheorem{lemma}[theo]{Lemma}

\newtheorem{coro}[theo]{Corollary}
\newtheorem{con}[theo]{Conjecture}
\newtheorem{prop}[theo]{Proposition}
\newtheorem{fact}[theo]{Fact}

\newcommand{\qed}{\hspace*{\fill} \rule{7pt}{7pt}}

\begin{document}
\maketitle
\begin{abstract}
An $r$-uniform graph $G$ is dense if and only if every proper subgraph $G'$ of $G$ satisfies $\lambda (G') < \lambda (G)$, where $\lambda (G)$ is the Lagrangian of a hypergraph $G$. In 1980's, Sidorenko  showed that $\pi(F)$,
the Tur\'an density of an $r$-uniform hypergraph $F$ is $r!$ multiplying the supremum of the Lagrangians of all dense $F$-hom-free $r$-uniform hypergraphs. This connection has been applied in estimating Tur\'an density of  hypergraphs. When $r=2$, the result of Motzkin and Straus shows that a  graph is dense if and only if it is a  complete graph. However, when $r\ge 3$, it becomes much harder to estimate the Lagrangians of $r$-uniform hypergraphs and to characterize the structure of all dense $r$-uniform graphs.
The main goal of this note is to give some sufficient conditions for $3$-uniform graphs with given substructures to be dense. For example, if $G$ is a $3$-graph with vertex set $[t]$ and $m$ edges containing $[t-1]^{(3)}$, then $G$ is dense if and only if $m \ge {t-1 \choose 3}+{t-2 \choose 2}+1$. We also give sufficient condition condition on the number of edges for a $3$-uniform hypergraph containing a large clique minus $1$ or $2$ edges to be dense.
\end{abstract}



\section{Introduction}
For a set $V$ and a positive integer $r$ we denote by $V^{(r)}$ the family of all $r$-subsets of $V$. An {\em $r$-uniform graph} or {\em $r$-graph $G$} consists of a set $V(G)$ of vertices and a set $E(G) \subseteq V(G) ^{(r)}$ of edges. An edge $e=\{a_1, a_2, \ldots, a_r\}$ will be simply denoted by $a_1a_2 \ldots a_r$. An $r$-graph $H$ is  a {\it subgraph} of an $r$-graph $G$, denoted by $H\subseteq G$ if $V(H)\subseteq V(G)$ and $E(H)\subseteq E(G)$. A subgraph $G$ {\em induced} by $V'\subseteq V$, denoted as $G[V']$, is the $r$-graph with vertex set $V'$ and edge set $E'=\{e\in E(G):e \subseteq V'\}$. Let ${\mathbb N}$ be the set of all positive integers. For $n \in  {\mathbb N}$, denote the set $\{1, 2, 3, \ldots, n\}$ by $[n]$. Let $K^{(r)}_t$ denote the {\em complete $r$-graph} on $t$ vertices, that is the $r$-graph on $t$ vertices containing all possible edges. A complete $r$-graph on $t$ vertices is also called {\em a clique} with order $t$. We also let $[n]^{(r)}$  denote the  complete $r$-uniform graph on the vertex set $[n]$. When $r=2$, an $r$-uniform graph is a {\em simple graph}.  When $r\ge 3$,  an $r$-graph is often called a {\em hypergraph}.
\bigskip

\begin{defi}
For  an $r$-uniform graph $G$ with the vertex set $[n]$,
edge set $E(G)$ and a weighting $\vec{x}=(x_1,\ldots,x_n) \in \mathbb R^n$,
define
$$\lambda (G,\vec{x})=\sum_{e \in E(G)}\prod\limits_{i\in e}x_{i}.$$
\end{defi}
The Lagrangian of
$G$, denoted by $\lambda (G)$, is defined as
 $$\lambda (G) = \max \{\lambda (G, \vec{x}): \vec{x} \in S \},$$
where $$S=\{\vec{x}=(x_1,x_2,\ldots ,x_n): \sum_{i=1}^{n} x_i =1, x_i
\ge 0 {\rm \ for \ } i=1,2,\ldots , n \}.$$

The value $x_i$ is called the {\em weight} of the vertex $i$ and any weighting $\vec{x} \in S$ is called a legal weighting.
A weighting $\vec{y}\in S$ is called an {\em optimum weighting} for $G$ if $\lambda (G, \vec{y})=\lambda(G)$.  The following fact is easily implied by the definition of the Lagrangian.

\begin{fact}\label{mono}
Let $G_1$, $G_2$ be $r$-uniform graphs and $G_1\subseteq G_2$. Then $\lambda (G_1) \le \lambda (G_2).$
\end{fact}

In \cite{MS}, Motzkin and Straus provided the following simple expression for the Lagrangian of a $2$-graph.

\begin{theo} {\rm (\cite{MS})} \label{MStheo}
If $G$ is a $2$-graph in which a largest clique has order $t$ then
$$\lambda(G)=\lambda(K^{(2)}_t)={1 \over 2}(1 - {1 \over t}).$$
\end{theo}

\begin{defi}
An $r$-uniform graph $G$ is dense if every proper subgraph $G'$ of $G$ satisfies $\lambda (G') < \lambda (G)$. This is equivalent to that all optimum weightings of $G$ are in the interior of $S$, in other words,  no coordinate in an optimum weighting is zero.
\end{defi}

In 1980's, Sidorenko \cite{Sidorenko} and Frankl and F\"uredi \cite{FF} developed the  Lagrangian method for hypergraph Tur\'an problems. Let $F$ and $G$ be hypergraphs. We say that a function $f:V(F)\rightarrow V(G)$ is a $homomorphism$ from hypergrah $F$ to hypergraph $G$ if it preserves edges, i.e. $f(i_1)f(i_2)\cdots f(i_k)\in E(G)$ for all $i_1i_2\cdots i_k\in E(F)$. We say that $G$ is $F$-$hom$-$free$ if there is no homomorphism from $F$ to $G$. Sidorenko \cite{Sidorenko} showed that $\pi(F)$,
the Tur\'an density of an $r$-uniform hypergraph $F$ is $r!$ multiplying the supremum of the Lagrangians of all dense $F$-hom-free $r$-uniform hypergraphs, i.e., $$\pi(F)=r! \sup \{\lambda(G): G \: {\rm is \: a \: dense} \: F\text{-}\hom\text{-} {\rm free} \:\: r\text{-} {\rm graph} \}.$$
Recent applications of this connection can be found in \cite{Keevash,HK,BIJ}.
When $r=2$, the result of Motzkin and Straus tells that all dense graphs are complete graphs. Therefore, for a graph $F$  with the chromatic number $t$, all dense $F$-hom-free graphs are complete graphs with order smaller than $t$. So $\pi(F)\le 2!\lambda(K_{t-1})=1-\frac{1}{t-1}$. This gives another proof of the fundamental result of Erd\H{o}s-Stone-Simonovits on Tur\'an densities of graphs. The key in the proof is that the result of Motzkin and Straus gives the structure of dense graphs and allows us to estimate the Lagrangian of a dense graph easily. However, when $r\ge 3$, it becomes much harder to estimate the Lagrangians of $r$-uniform hypergraphs and to characterize the structure of all dense $r$-uniform graphs. The main goal of this note is to characterize  $3$-graphs with certain given substructures to be dense.

In the following section, we give some elementary facts on dense $r$-graphs and some preliminary results needed in the proof. In sections 3, 4 and 5, we give results on $3$-uniform hypergraphs containing a large clique, a large clique with one edge or two edges removed correspondingly. Some open questions and remarks are given in section 6.

\section{Elementary Facts}
\begin{fact}{\rm (Cover pairs)} \label{fact1}
For every dense $r$-graph $G=(V,E)$ and for every pair $i,j\in V$, there exists an edge $e\in E(G)$ such that $\{i,j\} \subseteq e$.
\end{fact}
\noindent{\bf Proof.} Let $G=(V,E)$ be a dense $r$-graph with $n$ vertices. On the contrary suppose that there exists $i,j \in V$ such that $\{i,j\} \nsubseteq e$ for every $e\in E(G)$, i.e. $\{e \in E(G):\{i,j\}\subseteq e\}=\emptyset$. Let $\vec{x}=(x_{1},x_{2},\ldots ,x_{n})$ be an optimum weighting for $G$ with $k$ positive integers. Assume that $\frac{\partial \lambda (G, \vec{x})}{\partial x_i} \ge \frac{\partial \lambda (G, \vec{x})}{\partial x_j}$. Let $\vec{y}=(y_{1},y_{2},\ldots ,y_{n})$ be a weighting with $y_l=x_l$ for every $l\neq i,j$, $y_{i}=x_{i}+x_{j}$ and $y_{j}=0$. Then $\lambda(G,\vec{y})-\lambda(G,\vec{x})=x_{j}\left(\frac{\partial \lambda (G, \vec{x})}{\partial x_i}-\frac{\partial \lambda (G, \vec{x})}{\partial x_j}\right)-x_j^2\frac{\partial ^2 \lambda (G, \vec{x})}{\partial x_i \partial x_j}\ge 0$. It must be $\lambda(G,\vec{y})=\lambda(G,\vec{x})$ since $\vec{x}$ is an optimum weighting for $G$. This implies that $\lambda(G)=\lambda(G[V\setminus \{j\}])$, i.e. $G$ is not dense. This is a contradiction.
\qed

\begin{remark}\label{remark2.2}
Let $G$ be an $r$-graph on $t$ vertices. If $\lambda(G)> \lambda(K_{t-1}^{(r)})$ then $G$ is dense.
\end{remark}
\noindent{\bf Proof.}  $G$ is not dense if and only if there exists a subgraph $G'\subseteq K_{t-1}^{(r)}$ such that $\lambda(G)= \lambda(G')$. This implies that $\lambda(G) \le \lambda(K_{t-1}^{(r)})$, a contradiction.
\qed

\begin{fact}
Let $G$ be an $r$-graph on $t$ vertices with $m$ edges. If $m> \frac{t^r{t-1 \choose r}}{(t-1)^r}$ then $G$ is dense.
\end{fact}
\noindent{\bf Proof.} Let $G$ be an $r$-graph on $t$ vertices with $m> \frac{t^r{t-1 \choose r}}{(t-1)^r}$ edges. By Remark \ref{remark2.2}, it's sufficient to prove that $\lambda(G)>\lambda(K_{t-1}^{(r)})$. It's clear that $\lambda(K_{t-1}^{(r)})=\frac{{t-1 \choose r}}{(t-1)^r}$ since the weights of every vertex in the optimum weighting is $\frac{1}{t-1}$. Let ${\vec x}=(x_1, x_2, \ldots, x_t)$ be a legal weighting for $G$ with $x_j=\frac{1}{t}$ for every $j\in [t]$. Then
$\lambda(G)\ge \lambda(G,{\vec x})= \frac{m}{t^r}>\frac{{t-1 \choose r}}{(t-1)^r} =\lambda(K_{t-1}^{(r)})$.
\qed

\begin{remark}{\rm (Non-heredity)}
An induced subgraph of a dense $r$-graph may not be dense.
\end{remark}
\noindent{\bf Proof.} For example, $G=[t]^{(3)}\setminus \{(t-3)(t-1)t, (t-2)(t-1)t\}$ is dense by Theorem \ref{theo3}. In view of Fact \ref{fact1}, since $t-1$ and $t$ are not covered by any edge in $G[\{t-3,t-2,t-1,t\}]$, then $G[\{t-3,t-2,t-1,t\}]$ is not dense.
\qed
\bigskip

For an $r$-graph $G=(V,E)$ and $i\in V$, let $E_i=\{A \in V^{(r-1)}: A \cup \{i\} \in E\}$ denote the link of $i$. For a pair of vertices $i,j \in V$, let $E_{ij}=\{B \in V^{(r-2)}: B \cup \{i,j\} \in E\}$. Let $E^c_i=\{A \in V^{(r-1)}: A \cup \{i\} \in V^{(r)} \backslash E\}$ and
$E^c_{ij}=\{B \in V^{(r-2)}: B \cup \{i,j\} \in V^{(r)} \backslash E\}$. Denote $$E_{i\setminus j}=E_i\cap E^c_j, \ \ E_{(i_1i_2)\setminus (j_1j_2)}=E_{i_1i_2}\cap E^c_{j_1j_2}.$$

We say that an $r$-graph $G=(V,E)$  on vertex set $[n]$ is {\it left compressed} if $E_{j\setminus i}=\emptyset$ for any $1\le i<j\le n$. In other words, for any $i<j$, if $k_1 k_2\ldots k_{r-1}  \in E_j$, where $k_1, k_2, \ldots, k_{r-1} \neq i$, then  $k_1k_2 \ldots k_{r-1}  \in E_i$.

 An $r$-tuple $i_1 i_2\cdots i_r$ is called an {\it  ancestor} of an $r$-tuple $j_1j_2\cdots j_r$ if $i_1\le j_1$, $i_2\le j_2$, $\ldots$,  $i_r \le j_r$, and $i_1+i_2+ \cdots + i_r < j_1+j_2+ \cdots +j_r$. In this case, the  $r$-tuple $j_1j_2\cdots j_r$ is called a {\it descendant } of $i_1 i_2\cdots i_r$.  We say that $i_1 i_2\cdots i_r$ has higher hierarchy than $j_1j_2\cdots j_r$ if $i_1 i_2\cdots i_r$ is  an  ancestor of  $j_1j_2\cdots j_r$. We remark that an $r$-graph $G$ on the  vertex set $[n]$ is left-compressed if and only if all ancestors of an edge  are edges in $G$.
\bigskip

We will impose one additional condition on any optimum weighting ${\vec x}=(x_1, x_2, \ldots, x_n)$ for an $r$-graph $G$ in this paper:
\begin{eqnarray}
 &&|\{i : x_i > 0 \}|{\rm \ is \ minimal, i.e. \ if}  \ \vec y {\rm \ is \ a \ legal \ weighting \ for \ } G  {\rm \ satisfying }\nonumber \\
 &&|\{i : y_i > 0 \}| < |\{i : x_i > 0 \}|,  {\rm \  then \ } \lambda (G, {\vec y}) < \lambda(G) \label{conditionb}.
\end{eqnarray}

\bigskip
\begin{lemma} {\rm \cite{FR84}} \label{Lemma0} Let $G=(V,E)$ be an $r$-graph and ${\vec x}=(x_1, x_2, \ldots, x_n)$ be an optimum weighting for $G$ with $k$  positive weights $x_1, x_2, \ldots, x_k$. Then for every $\{i, j\} \in [k]^{(2)}$, {\rm(a)} $\lambda (E_i, {\vec x})=\lambda (E_j, \vec{x})=r\lambda(G)$, {\rm(b)} If ${\vec x}$ satisfies condition {\rm(\ref{conditionb})}, then there is an edge in $E$ containing both $i$ and $j$.
\end{lemma}

\begin{remark}\label{r1}  Let $G=(V,E)$ be an $r$-graph with the vertex set $[n]$ and ${\vec x}=(x_1, x_2, \ldots, x_n)$ be an optimum  weighting for $G$ with $k$  non-zero weights $x_1$, $x_2$, $\cdots$, $x_k$. Let $i, j$ be positive integers satisfying $1\le i<j\le k$. Then

{\rm(a)} Part {\rm(a)} in Lemma \ref{Lemma0} implies that
$$x_j\lambda(E_{ij}, {\vec x})+\lambda (E_{i\setminus j}, {\vec x})=x_i\lambda(E_{ij}, {\vec x})+\lambda (E_{j\setminus i}, {\vec x}).$$

{\rm(b)} If  $G$ is left-compressed, then
\begin{equation*}
x_i-x_j={\lambda (E_{i\setminus j}, {\vec x}) \over \lambda(E_{ij}, {\vec x})}
\end{equation*}
holds since $E_{j\setminus i}=\emptyset$. If  $G$ is left-compressed and  $E_{i\setminus j}=\emptyset$, then $x_i=x_j$.

{\rm(c)} If  $G$ is left-compressed, then
\begin{equation*}
x_1 \ge x_2 \ge \ldots \ge x_n \ge 0.
\end{equation*}
\end{remark}
\bigskip

\section{$3$-uniform hypergraphs containing a large clique}

The following result due to Peng-Zhao in \cite{PZ} implies that if $G$ is a $3$-graph with at most ${t-1 \choose 3} + {t-2 \choose 2}$ edges and $G$ contains a clique of order $t-1$, then $G$ is not dense.

\begin{theo}{\rm(\cite{PZ})} \label{theorempz} Let $m$ and $t$ be positive integers satisfying ${t-1 \choose 3} \le m \le {t-1 \choose 3} + {t-2 \choose 2}$. Let $G$ be a $3$-graph with $m$ edges and contain a clique of order  $t-1$. Then $\lambda(G) = \lambda([t-1]^{(3)})$.
\end{theo}

We give a sufficient condition for $r$-uniform hypergraphs containing a large clique to be dense.
\begin{theo}\label{theo3}
Let $m$, $t$ and $r\ge 3$ be positive integers satisfying $ m \ge {t-1 \choose r} + {t-2 \choose r-1}+1$. Let $G$ be an $r$-graph on $t$ vertices with $m$ edges. If $[t-1]^{(r)}\subseteq G$, then $G$ is dense.
\end{theo}
\noindent{\bf Proof.} Let $\vec{x}=(x_{1},x_{2},\ldots ,x_{t})$ be an optimum weighting of $G$. To complete the proof of the theorem, it's sufficient to prove that
$$\lambda(G)>\lambda(K_{t-1}^{(r)})={t-1 \choose r}\frac{1}{(t-1)^r}.$$
Consider a legal weighting $\vec{y}$ for $G$ with $y_{i}=\frac{1}{t-1}$ for every $i\in [t-2]$, $y_{t-1}=\frac{1-\delta}{t-1}$ and $y_{t}=\frac{\delta}{t-1}$, where $$0<\delta<\frac{1}{{t-2 \choose r-2}}.$$
Note that $\lambda(G,\vec{y})$ is the minimum if $\{A\cup (t-1)t:A\in {[t-2] \choose r-2}\}\subset E(G)$ when $\delta$ is small enough. Then
\begin{eqnarray*}
\lambda(G,\vec{y}) &\ge&{t-1 \choose r}\frac{1}{(t-1)^r}-{t-2 \choose r-1}\frac{\delta}{(t-1)^r}+{t-2 \choose r-2}\frac{(1-\delta)\delta}{(t-1)^r}\nonumber\\
&&+\left[{t-2 \choose r-1}+1-{t-2 \choose r-2}\right]\frac{\delta}{(t-1)^r}\nonumber\\
&=& {t-1 \choose r}\frac{1}{(t-1)^r}+\left[{t-2 \choose r-2}(1-\delta)-{t-2 \choose r-2}+1\right]\frac{\delta}{(t-1)^r}\nonumber\\
&=& {t-1 \choose r}\frac{1}{(t-1)^r}+\left[1-{t-2 \choose r-2}\delta\right]\frac{\delta}{(t-1)^r}\nonumber\\
&>& {t-1 \choose r}\frac{1}{(t-1)^r}\nonumber\\
&=& \lambda(K_{t-1}^{(r)}).
\end{eqnarray*}
By Remark \ref{remark2.2}, we have that $G$ is dense.
\qed

By Theorem \ref{theo3} and Theorem \ref{theorempz}, we get that
\begin{coro}\label{coro3}
Let $G$ be a $3$-graph with vertex set $[t]$ and $m$ edges containing $[t-1]^{(3)}$. Then $G$ is dense if and only if $m \ge {t-1 \choose 3}+{t-2 \choose 2}+1$.
\end{coro}
\bigskip

\section{$3$-uniform hypergraphs containing a large clique minus one edge}

Denote $H_1=[t-1]^{(3)} \setminus \{(t-3)(t-2)(t-1)\}$, $H_2=[t-1]^{(3)} \setminus \{(t-3)(t-2)(t-1), (t-4)(t-2)(t-1)\}$, $H_3=[t-1]^{(3)} \setminus \{(t-3)(t-2)(t-1), (t-5)(t-4)(t-1)\}$ and $H_4=[t-1]^{(3)} \setminus \{(t-3)(t-2)(t-1), (t-6)(t-5)(t-4)\}$. If a $3$-graph $G$ contains a clique of order $t-1$ minus one edge and $G$ does not contain a clique of order $t-1$, then $G[[t-1]]$ has only one non-isomorphic case, i.e. $G[[t-1]]=H_1$.
If a $3$-graph $G$ contains a clique of order $t-1$ minus two edges and $G$ does not contain a clique of order $t-1$ minus one edge then $G[[t-1]]$ has three non-isomorphic cases, i.e. $G[[t-1]]=H_2$, $H_3$ or $H_4$. Let $m$ and $t$ be positive integers satisfying $m \le {t-1 \choose 3} + {t-2 \choose 2}-2$, denote

$\lambda _{m,1} := \max\{\lambda(G):G$ is a $3$-graph with $m$ edges and $G[[t-1]]=H_1$ and $G$ doesn't contain a clique of order $t-1$\} and

$\lambda _{m,s} := \max\{\lambda(G):G$ is a $3$-graph with $m$ edges and $G[[t-1]]=H_s$  and $G$ doesn't contain a clique of order $t-1$ minus one edge\} for $s=2, 3 ,4$.

\subsection{Non-dense case}
In this subsection, we give a sufficient condition for $3$-uniform hypergraphs containing a large clique minus one edge to be non-dense.

\begin{theo}\label{theo11}
Let $m$ and $t\ge 5$ be positive integers satisfying ${t-1 \choose 3} \le m \le {t-1 \choose 3} + {t-2 \choose 2}-3$. Let $H$ be a $3$-graph with $m$ edges. If all but one edge of $[t-1]^{(3)}$ are in $H$  and $H$ doesn't contain a clique of order $t-1$, then $H$ is not dense.
\end{theo}
\noindent{\bf Proof.}  Without loss of generality, we can assume that $H[[t-1]]=H_1$.  We first prove the following lemma. The case of $s=2$ of Lemma \ref{Lemma10} will be applied to prove Theorem \ref{theo12}.
\begin{lemma}\label{Lemma10}
Let $m$ and $t$ be positive integers satisfying $m \le {t-1 \choose 3} + {t-2 \choose 2}-3s$, where $s=1,2$. Then there exists a left-compressed $3$-graph $G$ with $m$ edges satisfying $G[[t-1]]=H_s$ and $\lambda(G)= \lambda_{m,s}$.
\end{lemma}
\emph{Proof of Lemma \ref{Lemma10}.} Let $s=1$ or $2$. Let $F$ be a $3$-uniform graph with $m \le {t-1 \choose 3} + {t-2 \choose 2}-2s$ edges, $F[[t-1]]\cong H_s$ and $\lambda(F)= \lambda_{m,s}$. Let ${\vec x}=(x_1, x_2, \ldots, x_ n)$ be an optimum weighting of $F$ satisfying condition {\rm (1)} and $k$ be the number of non-zero weights in ${\vec x}$.
If $k\le t-1$, since the left-compressed $3$-graph $H_s \subseteq F$, then lemma holds. Now suppose that $k\ge t$. We can assume that $x_j\ge x_l$ when $1\le j<l\le t-1$ or $t\le j<l$, since otherwise we can just relabel the vertices of $F$ and obtain another extremal $3$-graph $F'$ satisfying $F'[[t-1]]\cong H_s$. Next we obtain a new $3$-graph $G$ from $F$ by performing the following:
\newline
1. If an edge $j_1j_2j_3$ has an ancestor  $i_1i_2i_3$ that is not in $F$, where $j_1,j_2,j_3<t$, then replace $j_1j_2j_3$ by $i_1i_2i_3$. Repeat this until there is no such an edge.
\newline
2. If an edge $j_1j_2t''$ in $F$ has an ancestor $i_1i_2t'$ that is not in $F$, where $t\le t'\le t''$, then replace $j_1j_2t''$ by $i_1i_2t'$. Repeat this until there is no such an edge.

Then $G$ satisfies the following properties:
\newline
1. The number of edges in $G$ is the same as the number of edges in $F$.
\newline
2. $\lambda(F)=\lambda(F,{\vec x})\le \lambda(G,{\vec x})\le \lambda(G)$.
\newline
3. $G[[t-1]]=H_s$.
\newline
4. All ancestors containing the vertex $t'$ of $j_1j_2t''$ in $G$ are in $G$, where $t\le t'\le t''$.

\noindent{Case 1. $s$=1.}
If $(t-3)(t-2)t \notin G$, then Properties 3 and 4 imply that  $G$ is left-compressed. Suppose that
$(t-3)(t-2)t \in G$.
 Note that $1(t-1)t\in G$ by Lemma \ref{Lemma0} (b) and $k\ge t$, we have
\begin{eqnarray*}
m &\ge& \left({t-1 \choose 3}-1\right)+{t-2 \choose 2}+1\nonumber\\
&=&{t-1 \choose 3}+{t-2 \choose 2},
\end{eqnarray*}
which contradicts to the number of edges of $G$.

\noindent{Case 2. $s$=2.} If $(t-4)(t-2)t \notin G$, then Properties 3 and 4 imply that  $G$ is left-compressed. Suppose that
$(t-4)(t-2)t \in G$. Then
\begin{eqnarray*}
m &\ge& \left({t-1 \choose 3}-2\right)+\left({t-2 \choose 2}-1\right)+1\nonumber\\
&=&{t-1 \choose 3}+{t-2 \choose 2}-2,
\end{eqnarray*}
which contradicts to the number of edges of $G$.
\qed

Now let us continue the proof of the theorem. By Lemma \ref{Lemma10} and Remark \ref{r1} (b), there exists a left-compressed $G$ such that  $G[[t-1]]=H_1$, $\lambda(G)= \lambda_{m,1}$ and  an optimum weighting $\vec{x}$ satisfying $x_1\ge x_2 \ge \cdots \ge x_k>x_{k+1}=\cdots=x_n=0$.
Since ${\vec x}$ has only $k$ positive weights, we can assume that $G$ is on $[k]$. If $k\le t-1$, then we are done. So suppose that $k\ge t$. We need the following lemma. The proof of the following lemma is similar to the proof of Lemma 2.5 in \cite{T} given by Talbot. For completeness we give a detailed proof later. For the cases of $s=2,3,4$ of Lemma \ref{Lemma2}, we need them to prove Theorem \ref{theo12}.
\begin{lemma}\label{Lemma2}
Let $l=1$ or $2$ and $s=1$, $2$, $3$ or $4$. Let $m_1$, $m_2$, $m_3$, $m_4$ and $t$ be positive integers satisfying $m_l \le {t-1 \choose 3} + {t-2 \choose 2}-2l$ and $m_2=m_3=m_4$. Let $G$ be a $3$-graph satisfying {\rm(a)} $G[[t-1]]=H_s$; {\rm(b)} if $i_1i_2t'$ is in $G$, then all ancestors of $i_1i_2t'$ containing the vertex $t''$ are in $G$, where $t\le t''\le t'$; and {\rm(c)} $\lambda(G)= \lambda_{m_s,s}$. Let $\vec{x}=(x_{1},x_{2},\ldots ,x_{n})$ be an optimum weighting of $G$ satisfying condition {\rm(1)} and $x_i\ge x_j$ if $i \le j$. Let $k$ be the number of positive weights in $\vec{x}$, then $$|[k-1]^{(3)}\backslash E| \le k+l-2.$$
\end{lemma}

Let us continue the proof of the theorem. Note that $G$ satisfies the condition of Lemma \ref{Lemma2}. By Lemma \ref{Lemma0} (b), $\vert E_{(k-1)k}\vert \ge 1$. If $k\ge t+1$, then applying Lemma \ref{Lemma2} ($l=1$), we have
\begin{eqnarray*}
m=\vert E\vert &=&\vert E\cap [k-1]^{(3)}\vert +\vert [k-2]^{(2)}\cap E_k \vert +\vert E_{(k-1)k}\vert \\
&\ge& {t \choose 3}-t+1 \\
&=& {t-1 \choose 3} + {t-2 \choose 2}-1,
\end{eqnarray*}
which contradicts that $m\le {t-1 \choose 3} + {t-2 \choose 2}-2$. Recall that $k\ge t$, so we have $k=t$, then $G$ is on $[t]$. To complete the proof, it's sufficient to prove the following Lemma. For the case of $s=2$ of Lemma \ref{Lemma13a}, we apply it to prove Theorem \ref{theo12}.

\begin{lemma}\label{Lemma13a}
Let $s=1$ or $2$. Let $m$, $t$ be positive integers satisfying $m_s \le {t-1 \choose 3} + {t-2 \choose 2}-3s$, and $G$ be a left-compressed $3$-graph on vertex set $[t]$ with $m$ edges. If $G[[t-1]]=H_s$, then $\lambda(G)=\lambda(H_s)$.
\end{lemma}
\emph{Proof of Lemma \ref{Lemma13a}.}
Let ${\vec x}=(x_1, x_2, \ldots, x_t)$ be an optimum weighting of $G$ and $k$ be the number of non-zero weights in ${\vec x}$.
Note that $x_1\ge x_2 \ge \cdots \ge x_t$.
Denote $b=|E_{(t-1)t}|$.
Since $G$ is left-compressed and $(t-3)(t-2)(t-1)\notin E(G)$, then $b \le t-4$.
It's sufficient to prove that $k\le t-1$. On the contrary suppose that $k=t$. By Remark \ref{r1} (b), we have
\begin{eqnarray}\label{eq1}
x_{t-2}-x_{t-1}&=&{\lambda (E_{(t-2) \setminus (t-1)}, {\vec x}) \over \lambda(E_{(t-2)(t-1)}, {\vec x})}\nonumber\\
&\le&\frac{x_t\lambda (E_{((t-2)t) \setminus ((t-1)t)}, {\vec x})}{\lambda(E_{(t-2)(t-1)}, {\vec x})}\nonumber\\
&\le& x_t
\end{eqnarray}
since if $ j \notin E_{(t-2)(t-1)}$, then $j \notin E_{(t-2)t}$, so $j \notin E_{((t-2)t) \setminus ((t-1)t)}$.
Similarly,
\begin{eqnarray}\label{eq2}
x_1-x_{t-1}&=&{\lambda (E_{1\setminus (t-1)}, {\vec x}) \over \lambda(E_{1(t-1)}, {\vec x})}\nonumber\\
&=&\frac{x_t\lambda (E_{(t-1)t}^{c}, {\vec x})+x_{t-2}\sum_{j=1}^{s}x_{t-j-2}}{1-x_{1}-x_{t-1}}\nonumber\\
&\le& x_{t}+\frac{s}{t-4}x_{t-2}.
\end{eqnarray}
Inequalities (\ref{eq1}) and (\ref{eq2}) imply that
\begin{eqnarray}\label{eq3}
x_1 \le \frac{t-4+s}{t-4}(x_{t-1}+x_t).
\end{eqnarray}
Again by Remark \ref{r1} (b) we have
\begin{eqnarray*}
x_{t-1}-x_{t}&=&{\lambda (E_{(t-1)\setminus t}, {\vec x}) \over \lambda(E_{(t-1)t}, {\vec x})}\nonumber\\
&=&\frac{\lambda (E_{t}^{c}\cap[t-2]^{(2)}, {\vec x})-\lambda (E_{t-1}^{c}\cap[t-2]^{(2)}, {\vec x})}{bx_{1}}.
\end{eqnarray*}
Since $$|[t-1]^{(3)}\cap E| ={t-1 \choose 3}-s$$
and
$$|E_{(t-1)t}|=b,$$
then
$$|E_{t}^{c}\cap[t-2]^{(2)}|={t-2 \choose 2}-\left(m-|[t-1]^{(3)}\cap E|-|E_{(t-1)t}|\right) \ge b+2s.$$
So
\begin{eqnarray}\label{eq4}
x_{t-1}-x_{t}&\ge& \frac{(b+s)x_{t-1}^2}{bx_{1}}\nonumber\\
&\ge& \frac{(b+s)x_{t-1}^2}{b \frac{t-4+s}{t-4}(x_{t-1}+x_t)}.
\end{eqnarray}
Inequality (\ref{eq4}) implies that $$\frac{b(t-4+s)}{t-4}(x_{t-1}^2-x_{t}^2)\ge (b+s)x_{t-1}^2.$$
So $$b(1+\frac{s}{t-4})>b+s.$$
This implies that
$$\frac{b}{t-4}>1,$$
so $$b>t-4,$$
which contradicts to $b\le t-4$.
\qed

Now let us continue the proof of the theorem. By Lemma \ref{Lemma10} and Lemma \ref{Lemma13a}, $\lambda(H) \le \lambda(G) = \lambda(H_1)= \lambda(H[[t-1]])$ and $H$ is not dense.
\qed

\bigskip

\subsection{Dense case}
In this subsection, we give a sufficient condition for $3$-uniform hypergraphs containing a large clique with one edge removed to be dense.

\begin{theo}\label{theo13}
Let $m$ and $t\ge 7$ be positive integers satisfying ${t-1 \choose 3} + {t-2 \choose 2}\le m \le {t \choose 3}$. Let $G$ be a $3$-graph with vertex set $[t]$ and $m$ edges. If all but one edge of $[t-1]^{(3)}$ are in $G$ and $G$ does not contain a clique of order $t-1$, then $G$ is dense.
\end{theo}
\noindent{\bf Proof.}
Let $G$ satisfy the condition of Theorem \ref{theo13}. Without loss of generality, we can assume that $G[[t-1]]=H_1$. First, we deduce the value of $\lambda(H_1)$. Denote $E=E(H_1)$. Note that $H_1$ is left-compressed. Let $\vec{x}=(x_{1},x_{2},\ldots ,x_{t-1})$ be an optimum weighting of $H_1$. By Remark \ref{r1}, we have $x_{1}=x_{2}=\ldots =x_{t-4}\ge x_{t-3}=x_{t-2}=x_{t-1}$. Assume that $x_{1}=x_{2}=\ldots =x_{t-4}=a$ and $x_{t-3}=x_{t-2}=x_{t-1}=b$, then $(t-4)a+3b=1$ and
\begin{eqnarray}\label{21}
\lambda(H_1)={t-4 \choose 3}a^3+3{t-4 \choose 2}a^2b+3(t-4)ab^2.
\end{eqnarray}
By Remark \ref{r1} (b) we have
\begin{eqnarray*}
a-b=x_1-x_{t-1}&=&{\lambda (E_{1\setminus (t-1)}, {\vec x}) \over \lambda(E_{1(t-1)}, {\vec x})}\nonumber\\
&=&\frac{x_{t-3}x_{t-2}}{1-x_1-x_{t-1}}\nonumber\\
&=&\frac{b^2}{1-a-b},
\end{eqnarray*}
then
\begin{eqnarray}\label{ab}
b=a-a^2.
\end{eqnarray}
Consider a legal weighting $\vec{y}$ for $G$ with $y_{i}=x_{i}=a$ for every $i\in [t-4]$, $y_{j}=x_j-\varepsilon=b-\varepsilon$ for $j=t-3,t-2,t-1$ and $y_t=3\varepsilon$, where $0< \varepsilon \ll b$. It's obvious that $y_{1}=y_{2}=\ldots =y_{t-4}\ge y_{t-3}=y_{t-2}=y_{t-1}\ge y_t$. Then
\begin{eqnarray*}
\lambda(G,\vec{y}) &\ge&{t-4 \choose 3}a^3+3{t-4 \choose 2}a^2(b-\varepsilon)+3(t-4)a(b-\varepsilon)^2+3(b-\varepsilon)^23\varepsilon+\nonumber\\
&&3(t-4)a(b-\varepsilon)\cdot3\varepsilon+\left[{t-2 \choose 2}+1-3(t-3)\right]a^2\cdot3\varepsilon\nonumber\\
&=& {t-4 \choose 3}a^3+3{t-4 \choose 2}a^2b+3(t-4)ab^2-3{t-4 \choose 2}a^2\varepsilon-6(t-4)ab\varepsilon\nonumber\\
&&+9b^2\varepsilon+9(t-4)ab\varepsilon+3\left[{t-2 \choose 2}+1-3(t-3)\right]a^2\varepsilon+\mathcal{O}(\varepsilon^{2})\nonumber\\
&=& \lambda(H_1)-3(t-3)a^2\varepsilon+3(t-4)ab\varepsilon+9b^2\varepsilon+\mathcal{O}(\varepsilon^{2}).
\end{eqnarray*}
Then we get that
\begin{eqnarray*}
&&\frac{1}{3}\left[\lambda(G,\vec{y})-\lambda(H_1)\right]\nonumber\\
&=&-(t-3)a^2\varepsilon+(t-4)ab\varepsilon+3b^2\varepsilon+\mathcal{O}(\varepsilon^{2})\nonumber\\
&=&-(t-3)a^2\varepsilon+(t-4)a(a-a^2)\varepsilon+3(a-a^2)^2\varepsilon+\mathcal{O}(\varepsilon^{2})\nonumber\\
&=& a^2\varepsilon \left[2-(t+2)a\right]+3a^4\varepsilon+\mathcal{O}(\varepsilon^{2}).
\end{eqnarray*}
Note that $(t-4)a+3b=1$ and $b=a-a^2$, then $3a^2-(t-1)a+1=0$. This implies that $a<0.2$ for every $t\ge 7$.
Hence
$2-(t+2)a=-3a^2-3a+1>0$ for every $a<0.2$. Hence $2-(t+2)a>0$ for every $t\ge 7$. So
\begin{eqnarray*}
&&\lambda(G,\vec{y})-\lambda(H_1)> 0
\end{eqnarray*}
when $\varepsilon$ is small enough. Note that every subgraph $G'$ of $G$ with order less than $t$ is a subgraph of $H_1$ under isomorphism,
then $\lambda(G')\le \lambda(H_1) < \lambda(G)$ and $G$ is dense.
\qed

\subsection{Special case}
Let $G$ be a $3$-graph on vertex set $[t]$ with $m$ edges  satisfying $G[[t-1]]=H_1$. Based on Theorems \ref{theo11} and \ref{theo13}, we can determine whether $G$ is dense according to the number of edges $m$ of $G$ except the case of $m={t-1 \choose 3} + {t-2 \choose 2}-1$ or $m={t-1 \choose 3} + {t-2 \choose 2}-2$.
For the case of $m={t-1 \choose 3} + {t-2 \choose 2}-2$, we believe that it would be non-dense, but we cannot prove it.
For the case of $m={t-1 \choose 3} + {t-2 \choose 2}-1$, a $3$-graph could be dense or non-dense. But we cannot give a characterization on a dense or non-dense $3$-graph in this case. We get the following result.
\begin{prop}\label{prop1}
Let $m$ and $t$ be positive integers satisfying $m={t-1 \choose 3} + {t-2 \choose 2}-1$. Let $G$ be a $3$-graph on $t$ vertices with $m$ edges satisfying $G[[t-1]]=H_1$. {\rm(a)} If $\{ijt:2\le i < j \le t-1 \}\subseteq E(G)$, then $G$ is not dense. {\rm(b)} If $|\{(t-3)(t-2)t,(t-3)(t-1)t,(t-2)(t-1)t\} \cap E(G)|\le 1$, then $G$ is dense.
\end{prop}
\noindent{\bf Proof.}
\noindent{(a)} Suppose $\{ijt:2\le i < j \le t-1 \}\subseteq E(G)$. Since $|\{ijt:2\le i < j \le t-1 \}|={t-2 \choose 2}$, $|E(G)\cap {[t-1] \choose 3}|={t-1 \choose 3}-1$ and $m={t-1 \choose 3} + {t-2 \choose 2}-1=|\{ijt:2\le i < j \le t-1 \}|+|E(G)\cap {[t-1] \choose 3}|$, then there is no edge else in $G$. So there is no edge $e$ of $G$ such that $\{1,t\}\subseteq e$. By Lemma \ref{Lemma0} (b), $G$ has an optimum weighting such that the number of its positive weights is less than $t$, i.e. $G$ is not dense.

\noindent{(b)} Suppose $|\{(t-3)(t-2)t,(t-3)(t-1)t,(t-2)(t-1)t\} \cap E(G)|\le 1.$  Let $\vec{x}=(x_{1},x_{2},\ldots ,x_{t-1})$ be an optimum weighting of $H_1$. By Remark \ref{r1} (c), we have $x_{1}=x_{2}=\ldots =x_{t-4}\ge x_{t-3}=x_{t-2}=x_{t-1}$. Assume that $x_{1}=x_{2}=\ldots =x_{t-4}=a$ and $x_{t-3}=x_{t-2}=x_{t-1}=b$. The same as the proof of Theorem \ref{theo13}, we have
$$(t-4)a+3b=1,$$
$$\lambda(H_1)={t-4 \choose 3}a^3+3{t-4 \choose 2}a^2b+3(t-4)ab^2$$
 and
$$ b=a-a^2. $$
Consider a legal weighting $\vec{y}$ for $G$ with $y_{i}=x_{i}=a$ for every $i\in [t-4]$, $y_{j}=x_j-\varepsilon=b-\varepsilon$ for $j=t-3,t-2,t-1$ and $y_t=3\varepsilon$, where $ \varepsilon $ is small enough. It's obvious that $y_{1}=y_{2}=\ldots =y_{t-4}\ge y_{t-3}=y_{t-2}=y_{t-1}\ge y_t$. Then
\begin{eqnarray*}
\lambda(G,\vec{y}) &\ge&{t-4 \choose 3}a^3+3{t-4 \choose 2}a^2(b-\varepsilon)+3(t-4)a(b-\varepsilon)^2+(b-\varepsilon)^23\varepsilon+\nonumber\\
&&3(t-4)a(b-\varepsilon)\cdot3\varepsilon+\left[{t-2 \choose 2}-3(t-4)-1\right]a^2\cdot3\varepsilon\nonumber\\
&=& \lambda(H_1)-3(t-4)a^2\varepsilon+3(t-4)ab\varepsilon+3b^2\varepsilon+\mathcal{O}(\varepsilon^{2}).
\end{eqnarray*}
Then we get that
\begin{eqnarray*}
&&\frac{1}{3}\left[\lambda(G,\vec{y})-\lambda(H_1)\right]\nonumber\\
&=&-(t-4)a^2\varepsilon+(t-4)ab\varepsilon+b^2\varepsilon+\mathcal{O}(\varepsilon^{2})\nonumber\\
&=&-(t-4)a^2\varepsilon+(t-4)a(a-a^2)\varepsilon+(a-a^2)^2\varepsilon+\mathcal{O}(\varepsilon^{2})\nonumber\\
\end{eqnarray*}
Recall that $b=a-a^2$ and $(t-4)a+3b=1$, so
\begin{eqnarray*}
&&-(t-4)a^2\varepsilon+(t-4)a(a-a^2)\varepsilon+(a-a^2)^2\varepsilon+\mathcal{O}(\varepsilon^{2})\nonumber\\
&=& a^2\varepsilon \left[1+a^2-(t-2)a\right]+\mathcal{O}(\varepsilon^{2})\nonumber\\
&=& a^2\varepsilon (3b+a^2-2a)+\mathcal{O}(\varepsilon^{2})\nonumber\\
&=& a^3(1-2a)\varepsilon+\mathcal{O}(\varepsilon^{2})\nonumber\\
&>& 0
\end{eqnarray*}
when $\varepsilon$ is small enough. Note that every subgraph $G'$ of $G$ with order less than $t$ is a subgraph of $H_1$ under isomorphism,
then $\lambda(G')\le \lambda(H_1) < \lambda(G)$ and $G$ is dense.
\qed

\subsection{The proof of Lemma \ref{Lemma2}}

\noindent{\bf Proof of Lemma \ref{Lemma2}.} Let $G=([n],E)$ and $\vec{x}$ satisfy the conditions. Denote $b=\max\{j: j(k-1)k\in E\}$. Since $G[[t-1]]=H_s$ and all ancestors containing the vertex $t''$ of $i_1i_2t'$ in $G$ are in $G$, where $t\le t''\le t'$ and $s=1,2,3,4$, we note that $b\le t-4$. Otherwise $(t-3)(t-2)t \in E(G)$, then $m\ge \left({t-1 \choose 3}-l\right)+{t-2 \choose 2}$, contradicts to $m \le {t-1 \choose 3} + {t-2 \choose 2}-2l$. Since $E_i=\{1, \ldots, i-1, i+1, \ldots, k\}^{(2)}$, for $1\le i\le b$, then $E_{i\setminus j}=\emptyset$ for $1\le i<j\le b$. Hence, by Remark \ref{r1}, we have $x_1=x_2=\cdots=x_b$. Also $x_i \ge x_j$ if $i<j$. Since $G[[t-1]]=H_s$, the conclusion is true if $k\le t$. So we may assume that $k\ge t+1$. We define a new legal weighting $\vec{y}$ for $G$ with $y_{j}=x_{j}$ for every $j\in [n]\setminus\{k-1,k\}$, $y_{k-1}=x_{k}+x_{k-1}$ and $y_{k}=0$. Note that $\lambda(E_{k-1},\vec{x})=\lambda(E_{k},\vec{x})$ from Lemma \ref{Lemma0} (a). So
\begin{eqnarray}\label{31}
\lambda(G,\vec{y})-\lambda(G,\vec{x})&=&x_{k}(\lambda(E_{k-1},\vec{x})-\lambda(E_{k},\vec{x}))-x_{k}^{2}\lambda(E_{(k-1)k},\vec{x})\nonumber\\
&=&-x_{k}^{2}\sum_{j=1}^{b}x_{j}\nonumber\\
&=&-bx_{1}x_{k}^{2}.
\end{eqnarray}
Since $y_{k}=0$ we may remove all edges containing the vertex $k$ from $E$ to form a new $3$-graph $G^{'}=([k-1],E^{'})$ with $\lambda(G',\vec{y})=\lambda(G,\vec{y}) $ and $E^{'}=E\cap[k-1]^{(3)}$. We will show that if Lemma \ref{Lemma2} fails to hold then there exists a set of edges $F\subset [k-1]^{(3)}\setminus E$ satisfying
\begin{eqnarray}\label{32}
\lambda(F,\vec{y})> bx_{1}x_{k}^{2},
\end{eqnarray}
\begin{eqnarray}\label{33}
|F|\le|E_k|
\end{eqnarray}
and there is no copy of a clique of order $t-1$ in $G'$ for $s=1$ and no copy of a clique of order $t-1$ minus one edge in $G[E'\cup F]$ for $s=2,3,4.$
Then, the graph $G^{''}=([k],E^{''})$, where $E^{''}=E^{'}\cup F$ satisfying $|E^{''}|\le |E|$ and
\begin{eqnarray*}
\lambda(G'',\vec{y})&=&\lambda(G',\vec{y})+\lambda(F,\vec{y})\nonumber\\
&>&\lambda(G',\vec{y})+ bx_{1}x_{k}^{2}\nonumber\\
&=&\lambda(G,\vec{x}).
\end{eqnarray*}
Hence $\lambda(G'')>\lambda(G)$. Note that $G''[[t-1]]=H_s$, which contradicts to $\lambda(G)=\lambda _{m_s,s}$.
We must construct the set of edges $F$ satisfying the above condition now. Since $E_{(k-1) \setminus 1}=\emptyset$, by Remark \ref{r1} we have $$x_1=x_{k-1}+{\lambda (E_{1\setminus (k-1)}, {\vec x}) \over \lambda(E_{1(k-1)}, {\vec x})}.$$
Hence
$$ bx_{1}x_{k}^{2}=bx_{k-1}x_{k}^{2}+\frac{bx_{k}^{3}\sum\limits_{j=b+1}^{k-2}x_{j}}{\sum\limits_{i=2,j\neq k-1}^{k}x_{j}}+ \frac{bx_{k}^{2}\lambda(C,\vec{x})}{\sum\limits_{j=2,j\neq k-1}^{k}x_{j}},$$
where $C=[k-2]^{(2)}\backslash E_{k-1}$.
Since $x_{j}\ge x_{l}$ when $j<l$ then
\begin{eqnarray}\label{35}
 bx_{1}x_{k}^{2}\le bx_{k-1}x_{k}^{2}\left[1+\frac{k-(b+2)}{k-3}\right]+\frac{bx_{k}\lambda(C,\vec{x})}{k-2}
\end{eqnarray}
Let $\alpha=\lceil\frac{b|C|}{k-2}\rceil$ and $\beta=b\lceil1+\frac{k-(b+2)}{k-3}\rceil$. Note that $\beta \le k-2$. Let $F'\subset [k-1]^{(3)}\setminus E$ consist of the first $\alpha$ heaviest edges in $[k-1]^{(3)}\setminus E$ containing the vertex $k-1$. Then
\begin{eqnarray*}
\lambda(F',\vec{y})\ge \frac{bx_{k}\lambda(C,\vec{x})}{k-2}+\alpha x_{k-1}x_{k}^{2}.
\end{eqnarray*}
Hence using (\ref{35})
\begin{eqnarray}\label{36}
\lambda(F',\vec{y})-bx_{1}x_{k}^{2}\ge (\alpha-\beta) x_{k-1}x_{k}^{2}.
\end{eqnarray}
We now distinguish two cases.

\noindent{\bf Case 1 $\alpha>\beta$.}

In this case $\lambda(F',\vec{y})>bx_{1}x_{k}^{2}$ so defining $F=F'$ satisfies (\ref{32}). We show that $|F|\le|E_k|$. Since $k\ge t$, $G[[t-1]]=H_s$ and all ancestors containing the vertex $t''$ of $i_1i_2t'$ in $G$ are in $G$, where $t\le t''\le t'$ and $s=1,2,3,4$, then $[b]^{(2)}\cup [b]\times\{b+1,\dots,k-1\}\subset E_k$. Hence
\begin{eqnarray}\label{37}
|E_k|\ge \frac{b[b-1+2(k-1-b)]}{2}\ge \frac{b(k-1)}{2}
\end{eqnarray}
since $b\le k-2$. Recall that $|F|=\alpha=\lceil\frac{b|C|}{k-2}\rceil$. Since $C\subset[k-2]^{(2)}$ we have $|C|\le {k-2 \choose 2}$. So using (\ref{37}) we obtain
$$|F|\le \left\lceil\frac{b(k-3)}{2}\right\rceil \le \frac{b(k-1)}{2} \le |E_k|.$$
Note that $G''[[t-1]]=H_s$. For $s=2,3,4$, if $G''$ contains a clique of order $t-1$ with $1$ edge removed ( this set of $t-1$ vertices would be different from $[t-1]$), then $|E(G'')|\ge {t-1 \choose 3}-2+{t-2 \choose 2}$, contradiction. Similarly, for $s=1$, $G''$ does not contain a clique of order $t-1$. So $F$ fulfills the purpose.

\noindent{\bf Case 2 $\alpha \le \beta$.}

Suppose Lemma \ref{Lemma2} fails, then $|[k-1]^{(3)}\setminus E|\ge k+l-1 \ge \beta+l+1$ (recall that $\beta\le k-2$). Let $F''\subset [k-1]^{(3)}\setminus E$ consist of any $\beta+1-\alpha$ edges in $[k-1]^{3}\setminus (E\cup F'\cup E(H_s)^c)$ and define $F=F'\cup F''$. Then since $\lambda(F'',\vec{y})\ge (\beta+1-\alpha)x_{k-1}^3$ and using (\ref{36}),
\begin{eqnarray*}
\lambda(F, \vec{y})-bx_{k-1}x_{k}^2&=&\lambda(F', \vec{y})-bx_{k-1}x_{k}^2+\lambda(F'', \vec{y})\nonumber\\
&\ge& (\beta+1-\alpha)x_{k-1}^3-x_{k-1}x_{k}^2(\beta-\alpha)\nonumber\\
&>&0.
\end{eqnarray*}
So (\ref{32}) is satisfied. We show that $\vert F\vert \le \vert E_k\vert.$ In fact,
$$\vert F\vert =\beta+1\le k-1\le {b(k-1) \over 2} \le \vert E_k\vert$$
when $b\ge 2$. If $b=1$, then
$$\vert F\vert =\beta+1=3\le k-2={b[b-1+2(k-1-b)] \over 2} \le \vert E_k\vert$$ since $k\ge t\ge 5$. As the same of Case 1, $F$ fulfills the purpose.
\qed

\section{$3$-uniform hypergraphs containing a large clique minus two edges}

\subsection{Non-dense case}
In this subsection, we give a sufficient condition for $3$-uniform hypergraphs containing a large clique minus two edges  to be  non-dense.

\begin{theo}\label{theo12}
Let $m$ and $t\ge 12$ be positive integers satisfying ${t-1 \choose 3} \le m \le {t-1 \choose 3} + {t-2 \choose 2}-6$. Let $H$ be a $3$-graph with $m$ edges. If all but two edges of $[t-1]^{(3)}$ are in $H$, then $H$ is not dense.
\end{theo}
\noindent{\bf Proof.}  Let $m$ and $t$ be positive integers satisfying $m \le {t-1 \choose 3} + {t-2 \choose 2}-6$. Let $H$ be a $3$-graph with $m$ edges satisfying all but two edges of $[t-1]^{(3)}$ are in $H$. Without loss of generality, we can assume that $H[[t-1]]=H_2$ or $H[[t-1]]=H_3$ or $H[[t-1]]=H_4$. It's sufficient to show that $\lambda(H)=\lambda(H[[t-1]])$.
We need the following lemma.
\begin{lemma}\label{Lemma101}
Let $s=3$ or $4$. Let $m$ and $t$ be positive integers satisfying $m \le {t-1 \choose 3} + {t-2 \choose 2}-6$. Then there exists a $3$-graph $G$ with $m$ edges satisfying
 {\rm (a)} $G[[t-1]]=H_s$,
 {\rm (b)} $G$ does not contain a subgraph isomorphic to $H_1$,
 {\rm (c)} if $i_1i_2t'$ is in $G$, then all ancestors of $i_1i_2t'$ containing the vertex $t''$ are in $G$, where $t\le t''\le t'$, and
 {\rm (d)} $\lambda(G)=\lambda_{m,s}$.
Moreover, there exists an optimum weighting $\vec{x}=(x_{1},x_{2},\ldots ,x_{n})$ of $G$ such that $x_{i}\ge x_{j}$ when $i<j$.
\end{lemma}
\emph{Proof of Lemma \ref{Lemma101}.} Let $s=3$ or $4$. Let $H$ be a $3$-uniform graph with $m \le {t-1 \choose 3} + {t-2 \choose 2}-6$ edges satisfying $H[[t-1]]\cong H_s$, $H$ does not contain a subgraph isomorphic to $H_1$, and $\lambda(H)= \lambda_{m,s}$. Let ${\vec x}=(x_1, x_2, \ldots, x_ n)$ be an optimum weighting of $H$ and $k$ be the number of non-zero weights in ${\vec x}$.
 We can assume that $x_1\ge x_2\ge \dots \ge x_{t-1}$ and $x_t \ge x_{t+1}\ge\dots \ge x_n$, otherwise we relabel the vertices.
 Next we obtain a new $3$-graph $G$ from $H$ by performing the following:
\newline
1. Case $s=3$. Note that there are exactly two elements $e,e' \in [t-1]^{(3)}\setminus E(H)$ satisfying $|e \cap e'|=1$. Replace  $e, e'$ by $(t-3)(t-2)(t-1), (t-5)(t-4)(t-1)$ respectively in $E(G)$.
\newline
Case $s=4$. Note that there are exactly two disjoint elements $e,e' \in [t-1]^{(3)}\setminus E(H)$. Replace  $e, e'$ by $(t-3)(t-2)(t-1), (t-6)(t-5)(t-4)$ respectively in $E(G)$.
\newline
2. If an edge $j_1j_2t''$ in $H$ has an ancestor $i_1i_2t'$  that is not in $E(H)$, where $t\le t'\le t''$, then replace $j_1j_2t''$ by $i_1i_2t'$. Repeat this until there is no such an edge.

Then $G$ satisfies the following properties:
\newline
1. The number of edges in $G$ is the same as the number of edges in $H$.
\newline
2. $\lambda(H)=\lambda(H,{\vec x})\le \lambda(G,{\vec x})\le \lambda(G)$.
\newline
3. $G[[t-1]]=H_s$.
\newline
4. All ancestors containing the vertex $t'$ of $i_1i_2t''$ in $G$ are in $G$, where $t\le t'\le t''$.
\newline
5. $G$ does not contain a subgraph isomorphic to $H_1$,

If $k\le t-1$, then taking the optimum vector of $G$ with the first $t-1$ coordinators the same as the coordinators of an optimum vector of $H_s$ and other coordinators being zero,  and the lemma holds. Now suppose that $k\ge t$.

 To complete the proof, it's sufficient to prove that $x_{t-1}\ge x_t$. By Remark \ref{r1} (a), we have $\lambda (E_{(t-1)\setminus t}, {\vec x})+x_t\lambda (E_{t(t-1)}, {\vec x})=\lambda (E_{t\setminus (t-1)}, {\vec x})+x_{t-1}\lambda (E_{t(t-1)}, {\vec x})$.
 So it's sufficient to prove that $\lambda (E_{(t-1)\setminus t}, {\vec x})-\lambda (E_{t\setminus (t-1)}, {\vec x})\ge 0$.

 For the case of $s=4$, note that $[t-1]^{(3)}\setminus E(G) = \{(t-3)(t-2)(t-1), (t-6)(t-5)(t-4)\}$ and by Property 4, we have $E_{t\setminus (t-1)} \subseteq \{(t-3)(t-2)\}$.
 Now we show that $E_{t\setminus (t-1)} = \emptyset$, i.e. $(t-3)(t-2) \notin E_{t\setminus (t-1)}$. Otherwise $(t-3)(t-2)t \in E(G)$, then $m \ge \left({t-1 \choose 3}-2\right)+ {t-2 \choose 2}$, which contradicts that  $m \le {t-1 \choose 3} + {t-2 \choose 2}-6$.
 So $\lambda (E_{(t-1)\setminus t}, {\vec x})-\lambda (E_{t\setminus (t-1)}, {\vec x})\ge 0$.

  For the case of $s=3$, note that $[t-1]^{(3)}\setminus E(G) = \{(t-3)(t-2)(t-1), (t-5)(t-4)(t-1)\}$ and by Property 4, we have $E_{t\setminus (t-1)} \subseteq \{(t-3)(t-2), (t-5)(t-4)\}$. Similar to  the case of $s=4$, we have  $(t-3)(t-2) \notin E_{t\setminus (t-1)}$.
  If $(t-5)(t-4)t \notin E(G)$, then $E_{t\setminus (t-1)}=\emptyset$ and we are done as above. So assume that $(t-5)(t-4)t \in E(G)$.
  Since $|E_{(t-1)\setminus t}|\ge |E_{(t-1)} \cap [t-2]^{(2)}|- |E_{t} \cap [t-2]^{(2)}|\ge {t-2 \choose 2}-2 -\left( m- ({t-1 \choose 3}-2)\right)\ge 2$,
  then $\{(t-4)(t-2), (t-5)(t-2)\} \subseteq E_{(t-1)\setminus t}$ or $\{(t-4)(t-2), (t-4)(t-3)\} \subseteq E_{(t-1)\setminus t}$.
  Hence $\lambda (E_{(t-1)\setminus t}, {\vec x})-\lambda (E_{t\setminus (t-1)}, {\vec x})\ge 2x_{t-4}x_{t-2}- x_{t-5}x_{t-4}$.
  To complete the proof, it's sufficient to prove that $ x_{t-5} \le 2x_{t-2}$.
Similar to the case of $s=4$, $(t-3)(t-1) \notin E_{t\setminus (t-2)}$, so $E_{t\setminus (t-2)}=\emptyset$. By Remark \ref{r1} (a) we have
\begin{eqnarray*}
x_{t-2}-x_{t}&=&{\lambda (E_{(t-2)\setminus t}, {\vec x})-\lambda (E_{t\setminus (t-2)}, {\vec x}) \over \lambda(E_{(t-2)t}, {\vec x})}\ge 0.
\end{eqnarray*}
Similarly,
\begin{eqnarray*}
x_{t-5}-x_{t-2}&=&{\lambda (E_{(t-5)\setminus (t-2)}, {\vec x})-\lambda (E_{(t-2)\setminus (t-5)}, {\vec x}) \over \lambda(E_{(t-5)(t-2)}, {\vec x})}.
\end{eqnarray*}
Since  $\lambda(H,\vec{x})= \lambda_{m,s}$, then  $\lambda(G,\vec{x})= \lambda_{m,s}$, that is, $\vec{x}$ is an optimal weighting of $G$.
If $x_{t+1} = 0$, then there is nothing to prove.
So assume that $x_{t+1} > 0$. We can also assume that $\vec{x}$ satisfying condition (1). By Lemma \ref{Lemma0} (b) we have $1t(t+1)\in E(G)$.
  So if $2(t-5)(t+1) \in E(G)$ (note that $(t-5)(t-4)t \in E(G)$), by Property 4 then $m \ge ({t-1 \choose 3}-2)+{t-5 \choose 2} + (t-1+t-7)>{t-1 \choose 3}+{t-5 \choose 2}-6$, a contradiction.
  So $2(t-5)(t+1) \notin E(G)$. 
  Note that $E_{(t-2)\setminus (t-5)}=\{(t-1)(t-4)\}$.
  Therefore
   $\lambda (E_{(t-5)\setminus (t-2)}, {\vec x})-\lambda (E_{(t-2)\setminus (t-5)}, {\vec x}) = x_{t-3}x_{t-1}+x_t \lambda (E_{(t-5)t}\cap E_{(t-2)t}^{c}, {\vec x})-x_{t-1}x_{t-4}$.
  Hence
 \begin{eqnarray*}
x_{t-5}-x_{t-2}=\frac{x_{t-3}x_{t-1}+x_t \lambda (E_{(t-5)t}\cap E_{(t-2)t}^{c}, {\vec x})-x_{t-6}x_{t-4}}{\lambda(E_{(t-5)(t-2)}, {\vec x})} \le x_t.
\end{eqnarray*}
Hence $x_{t-5} \le x_{t-2}+x_t \le 2x_{t-2}$.
\qed

\medskip
Let us continue the proof of the theorem. By Lemma \ref{Lemma10} and Lemma \ref{Lemma101}, we can assume that $H$ satisfies
(1) all ancestors containing the vertex $t'$ of $i_1i_2t''$ in $H$ are in $H$, where $t\le t'\le t''$
(2) $\lambda(H)= \lambda_{m,s}$ and
(3) $x_1\ge x_2 \ge \cdots \ge x_k>x_{k+1}=\cdots=x_n=0$. If $k\ge t+1$, then applying Lemma \ref{Lemma2}, we have
\begin{eqnarray*}
m=\vert E\vert &=&\vert E\cap [k-1]^{(3)}\vert +\vert [k-2]^{(2)}\cap E_k \vert +\vert E_{(k-1)k}\vert \\
&\ge & {t \choose 3}-t = {t-1 \choose 3} + {t-2 \choose 2}-2,
\end{eqnarray*}
which contradicts the assumption that $m\le {l-1 \choose 3} + {l-2 \choose 2}-6$. Recall that $k\ge t$, so we have $k=t,$ then $H$ is on $[t]$. We prove the following lemma.

\begin{lemma}\label{Lemma13b}
Let $s=3$ or $4$, $m$ and $t\ge 12$ be positive integers satisfying $m \le {t-1 \choose 3} + {t-2 \choose 2}-6$. Let $G$ be a $3$-graph on vertex set $[t]$ with $m$ edges satisfying $G[[t-1]]=H_s$ and all ancestors containing the vertex $t$ of $i_1i_2t$ in $G$ are in $G$. And $G$ has an optimum weighting $\vec{x}=(x_{1},x_{2},\ldots ,x_{n})$ satisfying $x_{i}\ge x_{j}$ when $i<j$. Then $\lambda(G)=\lambda(G[[t-1]])$.
\end{lemma}
\emph{Proof of Lemma \ref{Lemma13b}.}
\noindent{\bf Case $s=3$}.
Note that $(t-3)(t-2)t\notin E(G)$, otherwise $m \ge {t-1 \choose 3}-2 + {t-2 \choose 2} > {t-1 \choose 3} + {t-2 \choose 2}-6$, a contradiction.  Let $k$ be the number of non-zero weights in ${\vec x}$. Denote $b=|E_{(t-1)t}|$, clearly $b\le t-2$.
If $(t-5)(t-4)t\in E(G)$, then $i(t-1)t$, $i(t-2)t$, $i(t-3)t$ and $i'j't$, where $1 \le i \le b$ and $1 \le i' < j' \le t-4$, are edges in $G$. Therefore
 $3b+{t-4 \choose 2}\le m- \left({t-1 \choose 3}-2\right) \le {t-2 \choose 2}-4$,
 so $b\le \frac{2}{3}t-\frac{11}{3} \le \frac{2}{3}t-3$.
It's sufficient to prove that $k\le t-1$. On the contrary suppose that $k=t$. Since $E_{(t-1)\setminus (t-5)}=\emptyset$, by Remark \ref{r1} (a) we have
\begin{eqnarray*}
x_{t-5}-x_{t-1}&=&{\lambda (E_{(t-5)\setminus (t-1)}, {\vec x}) \over \lambda(E_{(t-5)(t-1)}, {\vec x})}\nonumber\\
&=&\frac{x_t\lambda (E_{(t-5)t\setminus(t-1)t}, {\vec x})+x_{t-3}x_{t-2}}{\lambda(E_{(t-5)(t-1)}, {\vec x})}.
\end{eqnarray*}
If $ j \notin E_{(t-5)(t-1)}$, then $j \notin E_{(t-5)t \setminus (t-1)t}$.
Note that $\lambda(E_{(t-5)(t-1)}, {\vec x}) \ge 1- x_{t-5}-x_{t-4}-x_{t-1}-x_t \ge (t-4)x_{t-2}$.
Then
\begin{eqnarray*}
 x_{t-5}-x_{t-1} \le  x_{t}+\frac{x_{t-3}}{t-4}.
\end{eqnarray*}
Note that $ x_{t-3}\le x_{t-5}$, so
\begin{eqnarray}\label{eq12}
x_{t-5} \le \frac{t-4}{t-5}(x_{t-1}+x_t).
\end{eqnarray}
Similarly,
\begin{eqnarray}\label{eq13}
x_1-x_{t-1}&=&{\lambda (E_{1\setminus (t-1)}, {\vec x}) \over \lambda(E_{1(t-1)}, {\vec x})}\nonumber\\
&=&\frac{x_t\lambda (E_{1t\setminus(t-1)t}, {\vec x})+x_{t-3}x_{t-2}+x_{t-5}x_{t-4}}{1-x_1-x_{t-1}}\nonumber\\
&\le& x_{t}+\frac{2x_{t-5}}{t-5}.
\end{eqnarray}
Inequalities (\ref{eq12}) and (\ref{eq13}) imply that
\begin{eqnarray*}\label{eq3}
x_1 \le \left(1+\frac{2(t-4)}{(t-5)^2}\right)(x_{t-1}+x_t).
\end{eqnarray*}
If $E_{t\setminus (t-1)}=\emptyset$, then  $b\le t-6$.
By Remark \ref{r1} (a) we have
\begin{eqnarray*}
x_{t-1}-x_{t}&=&{\lambda (E_{(t-1)\setminus t}, {\vec x}) \over \lambda(E_{(t-1)t}, {\vec x})}\nonumber\\
&=&\frac{\lambda (E_{t}^{c}\cap[t-2]^{(2)}, {\vec x})-\lambda (E_{t-1}^{c}\cap[t-2]^{(2)}, {\vec x})}{bx_{1}}.
\end{eqnarray*}
Since $$|[t-1]^{(3)}\cap E| ={t-1 \choose 3}-2$$
and
$$|E_{(t-1)t}|=b,$$
then
$$|E_{t}^{c}\cap[t-2]^{(2)}|={t-2 \choose 2}-\left(m-|[t-1]^{(3)}\cap E|-|E_{(t-1)t}|\right) \ge b+4.$$
So
\begin{eqnarray}\label{eq4}
x_{t-1}-x_{t}&\ge& \frac{(b+2)x_{t-1}^2}{bx_{1}}\nonumber\\
&\ge& \frac{(b+2)x_{t-1}^2}{b \left(1+\frac{2(t-4)}{(t-5)^2}\right)(x_{t-1}+x_t)}.
\end{eqnarray}
Inequality (\ref{eq4}) implies that $$b \left(1+\frac{2(t-4)}{(t-5)^2}\right)(x_{t-1}^2-x_{t}^2)\ge (b+2)x_{t-1}^2.$$
This implies that
$b\ge t-5$ which contradicts to $b\le t-6$.

So we assume that $E_{t\setminus (t-1)}=\{(t-5)(t-4)\}$. Again by Remark \ref{r1} (a) we have
\begin{eqnarray*}
x_{t-1}-x_{t}&=&{\lambda (E_{(t-1)\setminus t}, {\vec x})-\lambda (E_{t\setminus (t-1)}, {\vec x}) \over \lambda(E_{(t-1)t}, {\vec x})}\nonumber\\
&=&\frac{\lambda (E_{t}^{c}\cap[t-2]^{(2)}, {\vec x})-x_{t-3}x_{t-2}-x_{t-5}x_{t-4}}{bx_{1}}.
\end{eqnarray*}
Since there is at least one element containing $t-4$ in $E_{t}^{c}\cap[t-2]^{(2)}$, $(t-3)(t-2) \in E_{t}^{c}\cap[t-2]^{(2)}$ and
$$|[t-1]^{(3)}\cap E| ={t-1 \choose 3}-2, \:\: |E_{(t-1)t}|=b,$$
then
$$|E_{t}^{c}\cap[t-2]^{(2)}|={t-2 \choose 2}-\left(m-|[t-1]^{(3)}\cap E|-|E_{(t-1)t}|\right) \ge b+4.$$
So
\begin{eqnarray*}
x_{t-1}-x_{t}&\ge& \frac{(b+2)x_{t-2}x_{t-1}-(x_{t-5}-x_{t-1})x_{t-4}}{bx_{1}}\nonumber\\
&\ge&\frac{(b+2)x_{t-2}x_{t-1}-\frac{t-4}{t-5}(\frac{1}{t-5}x_{t-1}+\frac{t-4}{t-5}x_{t})(x_{t-1}+x_{t})}{b(1+\frac{2(t-4)}{(t-5)^2})(x_{t-1}+x_t)}.
\end{eqnarray*}
Multiplying both sides by the denominator, applying $x_{t-2}\ge x_{t-1}$ and combining common terms, we get
\begin{eqnarray}\label{eq14}
\left(\frac{(2b+1)(t-4)}{(t-5)^2}-2\right)x_{t-1}^2+\frac{(t-4)(t-3)}{(t-5)^2}x_{t-1}x_{t}+\left( \frac{(t-4)^2-2b(t-4)}{(t-5)^2}-b\right)x_{t}^2 \ge 0.
\end{eqnarray}
Now we consider two cases according to whether $x_{t}\ge \frac{1}{2}x_{t-1}$.

Case 1. $x_{t}\ge \frac{1}{2}x_{t-1}$.
Since $x_{t}\le x_{t-1}$, then (\ref{eq14}) implies that
$$\left(\frac{(2b+1)(t-4)}{(t-5)^2}-2\right)x_{t-1}^2+\frac{(t-4)(t-3)}{(t-5)^2}x_{t-1}^2+\left( \frac{(t-4)^2-2b(t-4)}{(t-5)^2}-b\right)x_{t}^2 \ge 0.$$
For every $t\ge 12$ and $1\le b \le \frac{2}{3}t-3$, observing that $\frac{(t-4)^2-2b(t-4)}{(t-5)^2}-b < 0$. Applying  $x_{t}\ge \frac{1}{2}x_{t-1}$ to the above inequality we obtained that
$$ \left(\frac{(t-4)(t+2b-2)}{(t-5)^2}-2+\frac{(t-4)^2-b(t-5)^2-2b(t-4)}{4(t-5)^2}\right)x_{t-1}^2 \ge 0.$$
The left of the above inequality is negative. This is a contradiction.

Case 2. $x_{t}< \frac{1}{2}x_{t-1}$.
Then (\ref{eq14}) implies that
$$ \left(\frac{(2b+1)(t-4)+\frac{1}{2}(t-4)(t-3)}{(t-5)^2}-2\right)x_{t-1}^2> \left(b+ \frac{2b(t-4)-(t-4)^2}{(t-5)^2}\right)x_{t}^2.$$
Note that $1 \le b\le \frac{2}{3}t-3$ and $t\ge 12$, the left of the above inequality is negative and the right of the above inequality is positive. This is a contradiction.

\noindent{\bf Case $s=4$.} (This proof is very similar to the proof of Lemma \ref{Lemma13a})
Let $G$ be a $3$-graph on vertex set $[t]$ with $m$ edges satisfying $H[[t-1]]=H_4$ and all ancestors containing the vertex $t$ of $i_1i_2t$ in $G$ are in $G$. Let ${\vec x}=(x_1, x_2, \ldots, x_t)$ be an optimum weighting of $G$ and $k$ be the number of non-zero weights in ${\vec x}$. Denote $b=|E_{(t-1)t}|$. Note that $x_1\ge x_2 \ge \cdots \ge x_t$, $b\le t-3$ and $(t-3)(t-2)t\notin E(G)$. It's sufficient to prove that $k\le t-1$. On the contrary suppose that $k=t$. Note that $E_{(t-1) \setminus (t-2)}=\emptyset$, by Remark \ref{r1} (a) we have
\begin{eqnarray}\label{eq31}
x_{t-2}-x_{t-1}&=&{\lambda (E_{(t-2) \setminus (t-1)}, {\vec x}) \over \lambda(E_{(t-2)(t-1)}, {\vec x})}\nonumber\\
&\le&\frac{x_t\lambda (E_{(t-2)t \setminus (t-1)t}, {\vec x})}{\lambda(E_{(t-2)(t-1)}, {\vec x})}\nonumber\\
&\le& x_t
\end{eqnarray}
since if $ j \notin E_{(t-2)(t-1)}$, then $j =t$ or $t-3$, so $j \notin E_{(t-2)t \setminus (t-1)t}$. Since $E_{1\setminus(t-1)}=\emptyset$, then
\begin{eqnarray}\label{eq32}
x_1-x_{t-1}&=&{\lambda (E_{1\setminus (t-1)}, {\vec x}) \over \lambda(E_{1(t-1)}, {\vec x})}\nonumber\\
&=&\frac{x_t\lambda (E_{(t-1)t}^{c}, {\vec x})+x_{t-3}x_{t-2}}{1-x_{1}-x_{t-1}}\nonumber\\
&\le& x_{t}+\frac{1}{t-4}x_{t-2}.
\end{eqnarray}
Inequalities (\ref{eq31}) and (\ref{eq32}) imply that
\begin{eqnarray*}\label{eq33}
x_1 \le \frac{t-3}{t-4}(x_{t-1}+x_t).
\end{eqnarray*}
Note that $(t-3)(t-2)t\notin E(G)$, by Remark \ref{r1} (a) again we have
\begin{eqnarray*}
x_{t-1}-x_{t}&=&{\lambda (E_{(t-1)\setminus t}, {\vec x}) \over \lambda(E_{(t-1)t}, {\vec x})}\nonumber\\
&=&\frac{\lambda (E_{t}^{c}\cap[t-2]^{(2)}, {\vec x})-x_{t-3}x_{t-2}}{bx_{1}}.
\end{eqnarray*}
Since $$|[t-1]^{(3)}\cap E| ={t-1 \choose 3}-2$$
and
$$|E_{(t-1)t}|=b,$$
then
$$|E_{t}^{c}\cap[t-2]^{(2)}|={t-2 \choose 2}-\left(m-|[t-1]^{(3)}\cap E|-|E_{(t-1)t}|\right) \ge b+4.$$
So
\begin{eqnarray}\label{eq34}
x_{t-1}-x_{t}&\ge& \frac{(b+3)x_{t-1}^2}{bx_{1}}\nonumber\\
&\ge& \frac{(b+3)x_{t-1}^2}{b \frac{t-3}{t-4}(x_{t-1}+x_t)}.
\end{eqnarray}
Inequality (\ref{eq34}) implies that $$\frac{b(t-3)}{t-4}(x_{t-1}^2-x_{t}^2)\ge (b+3)x_{t-1}^2.$$
So $$b(1+\frac{1}{t-4})>b+3.$$
This implies that
$$b>3(t-4),$$
which contradicts to $b\le t-2$. \qed

Then by Lemma \ref{Lemma13a} and Lemma \ref{Lemma13b}, $\lambda(H) \le \lambda(G)\le \lambda(H_1)= \lambda(H[[t-1]])$ and $H$ is not dense. This completes the proof.
\qed

\subsection{Dense case}
In this subsection, we give a sufficient condition for $3$-uniform hypergraphs containing a large clique minus two edges to be  dense.

\begin{theo}\label{theo14}
Let $m$ and $t\ge 9$ be positive integers satisfying ${t-1 \choose 3} + {t-2 \choose 2}-1\le m \le {t \choose 3}$. Let $G$ be a $3$-graph with vertex set $[t]$ and $m$ edges. If all but two edges of $[t-1]^{(3)}$ are in $G$ and $G$ does not contain a clique of order $t-1$ minus one edge, then $G$ is dense.
\end{theo}
\noindent{\bf Proof.}
Let $A,B$ be the two edges satisfying $A,B\in [t-1]^{(3)}\setminus E$. There are three cases, $|A\cap B|=0$, $|A\cap B|=1$ and $|A\cap B|=2$. Let $G_1,G_2,G_3=[t-1]^{(3)}\setminus \{A,B\}$ be the $3$-graphs satisfying $|A\cap B|=0$, $|A\cap B|=1$ and $|A\cap B|=2$ accordingly. Without loss of generality, we can assume that $G_1=H_4$, $G_2=H_3$ and $G_3=H_2$.

\noindent{Case 1. $|A\cap B|=0$.}
First, we deduce the value of $\lambda (H_4)$. Let $\vec{x}=(x_{1},x_{2},\ldots ,x_{t-1})$ be an optimum weighting of $H_4$. By Remark \ref{r1} we have $x_{1}=x_{2}=\ldots =x_{t-7}\ge x_{t-6}=\ldots =x_{t-1}$. Assume that $x_{1}=x_{2}=\ldots =x_{t-7}=a$ and $x_{t-6}=\ldots =x_{t-1}=b$, then $(t-7)a+6b=1$ and \begin{eqnarray}\label{22}
\lambda(H_4)={t-7 \choose 3}a^3+6{t-7 \choose 2}a^2b+(t-7){6 \choose 2}ab^2+18b^3.
\end{eqnarray}
By Remark \ref{r1} (a) we have
\begin{eqnarray*}
a-b=x_1-x_{t-1}&=&{\lambda (E_{1\setminus (t-1)}, {\vec x}) \over \lambda(E_{1(t-1)}, {\vec x})}\nonumber\\
&=&\frac{x_{t-3}x_{t-2}}{1-x_1-x_{t-1}}\nonumber\\
&=&\frac{b^2}{1-a-b}\nonumber\\
&<&\frac{b}{t-3},
\end{eqnarray*}
where $E$ is the edge set of $H_4$. So
\begin{eqnarray}\label{ab1}
b>\frac{t-3}{t-2}a.
\end{eqnarray}
Consider a legal weighting $\vec{y}=(y_{1},y_{2},\ldots ,y_{t})$ for $G$ with $y_{i}=x_{i}$ for every $i\neq 1$, $y_{1}=x_1-\varepsilon$ and $y_{t}=\varepsilon$, where $0< \varepsilon <a$. Since $a>b$ and $\lambda(G,\vec{y})$ is the minimum when the number of edges containing the vertex $t$ with weight $a^2\varepsilon$ is minimum. Then we get that
\begin{eqnarray*}
&&\lambda(G,\vec{y})-\lambda(H_4)\nonumber\\
&\ge& -{t-8 \choose 2}a^2\varepsilon-6(t-8)ab\varepsilon-{6 \choose 2}b^2\varepsilon+{6 \choose 2}b^2\varepsilon+6(t-7)ab\varepsilon \nonumber\\
&&+\left[{t-2 \choose 2}+1-6(t-7)-15\right]a^2\varepsilon+\mathcal{O}(\varepsilon^{2})\nonumber\\
&=&-5a^2\varepsilon+6ab\varepsilon+\mathcal{O}(\varepsilon^{2})\nonumber\\
&>& \frac{(t-8)a^2\varepsilon}{t-2}+\mathcal{O}(\varepsilon^{2})\nonumber\\
&>&0,
\end{eqnarray*}
when $\varepsilon$ is small enough.

\noindent{Case 2. $|A\cap B|=1$.}
First, we deduce the value of $\lambda(H_3)$. Let $\vec{x}=(x_{1},x_{2},\ldots ,x_{t-1})$ be an optimum weighting of $H_3$. By Remark \ref{r1} we have $x_{1}=x_{2}=\ldots =x_{t-6}\ge x_{t-5}=x_{t-4}=x_{t-3}= x_{t-2}\ge x_{t-1}$. Assume that $x_{1}=x_{2}=\ldots =x_{t-6}=a$, $x_{t-5}=x_{t-4}=x_{t-3}= x_{t-2}=b$ and $x_{t-1}=c$, then $(t-6)a+4b+c=1$ and
\begin{eqnarray*}\label{42}
\lambda(H_3)={t-6 \choose 3}a^3+4{t-6 \choose 2}a^2b+{t-6 \choose 2}a^2c+ 6(t-6)ab^2+4(t-6)abc+4b^2c.
\end{eqnarray*}
By Remark \ref{r1} (a) we have
\begin{eqnarray*}
a-b=x_1-x_{t-2}&=&{\lambda (E_{1\setminus (t-2)}, {\vec x}) \over \lambda(E_{1(t-2)}, {\vec x})}\nonumber\\
&=&\frac{x_{t-3}x_{t-1}}{1-x_1-x_{t-2}}\nonumber\\
&=&\frac{bc}{1-a-b}\nonumber\\
&<&\frac{b}{t-3},
\end{eqnarray*}
where $E$ is the edge set of $H_3$. So
\begin{eqnarray}\label{ab3}
b>\frac{t-3}{t-2}a.
\end{eqnarray}
And
\begin{eqnarray*}
b-c=x_{t-2}-x_{t-1}&=&{\lambda (E_{(t-2)\setminus (t-1)}, {\vec x}) \over \lambda(E_{(t-2)(t-1)}, {\vec x})}\nonumber\\
&=&\frac{x_{t-5}x_{t-4}}{1-x_{t-3}-x_{t-2}-x_{t-1}}\nonumber\\
&=&\frac{b^2}{1-2b-c}\nonumber\\
&<&\frac{b}{t-4},
\end{eqnarray*}
then together with (\ref{ab3}) we have
\begin{eqnarray}\label{ac1}
c>\frac{t-5}{t-4}b>\frac{t-5}{t-2}a.
\end{eqnarray}
Consider a legal weighting $\vec{y}=(y_{1},y_{2},\ldots ,y_{t})$ for $G$ with $y_{i}=x_{i}$ for every $i\neq 1$, $y_{1}=x_1-\varepsilon$ and $y_{t}=\varepsilon$, where $0< \varepsilon <a$. Since $a>b>c$ and $\lambda(G)$ is the minimum when the number of edges containing the vertex $t$ with weight $a^2\varepsilon$ is minimum. Then we get that
\begin{eqnarray*}
&&\lambda(G,\vec{y})-\lambda(H_3)\nonumber\\
&\ge& -{t-7 \choose 2}a^2\varepsilon-4(t-7)ab\varepsilon-(t-7)ac\varepsilon-6b^2\varepsilon-4bc\varepsilon+4bc\varepsilon+(t-6)ac\varepsilon\nonumber\\
&&+6b^2\varepsilon+4(t-6)ab\varepsilon +\left[{t-2 \choose 2}+1-5(t-6)-10\right]a^2\varepsilon+\mathcal{O}(\varepsilon^{2})\nonumber\\
&=&-4a^2\varepsilon+4ab\varepsilon+ac\varepsilon+\mathcal{O}(\varepsilon^{2})\nonumber\\
&>& \frac{(t-9)a^2\varepsilon}{t-2}+\mathcal{O}(\varepsilon^{2})\nonumber\\
&>&0
\end{eqnarray*}
when $\varepsilon$ is small enough.

\noindent{Case 3. $|A\cap B|=2$.}
First, we deduce the value of $\lambda(H_2)$. Let $\vec{x}=(x_{1},x_{2},\ldots ,x_{t-1})$ be an optimum weighting of $H_2$. By Remark \ref{r1} we have $x_{1}=x_{2}=\ldots =x_{t-5}\ge x_{t-4}=x_{t-3}\ge x_{t-2}=x_{t-1}$. Assume that $x_{1}=x_{2}=\ldots =x_{t-5}=a$, $x_{t-4}=x_{t-3}=b$ and $x_{t-2}=x_{t-1}=c$, then $(t-5)a+2b+2c=1$ and
\begin{eqnarray*}\label{22}
\lambda(H_2)={t-5 \choose 3}a^3+2{t-5 \choose 2}a^2b+2{t-5 \choose 2}a^2c+(t-5)ab^2+4(t-5)abc+(t-5)ac^2+2b^2c.
\end{eqnarray*}
By Remark \ref{r1} (a) we have
\begin{eqnarray*}
a-b=x_1-x_{t-4}&=&{\lambda (E_{1\setminus (t-4)}, {\vec x}) \over \lambda(E_{1(t-4)}, {\vec x})}\nonumber\\
&=&\frac{x_{t-2}x_{t-1}}{1-x_1-x_{t-4}}\nonumber\\
&=&\frac{c^2}{1-a-b}\nonumber\\
&<&\frac{c}{t-3}\nonumber\\
&<&\frac{b}{t-3},
\end{eqnarray*}
where $E$ is the edge set of $H_2$. So
\begin{eqnarray}\label{ab1}
b>\frac{t-3}{t-2}a.
\end{eqnarray}
And
\begin{eqnarray*}
a-c=x_1-x_{t-1}&=&{\lambda (E_{1\setminus (t-1)}, {\vec x}) \over \lambda(E_{1(t-1)}, {\vec x})}\nonumber\\
&=&\frac{(x_{t-4}+x_{t-3})x_{t-2}}{1-x_1-x_{t-1}}\nonumber\\
&=&\frac{2bc}{1-a-c}\nonumber\\
&<&\frac{2c}{t-4},
\end{eqnarray*}
then
\begin{eqnarray}\label{ac}
c>\frac{t-4}{t-2}a.
\end{eqnarray}
Consider a legal weighting $\vec{y}=(y_{1},y_{2},\ldots ,y_{t})$ for $G$ with $y_{i}=x_{i}$ for every $i\in [t-3]$, $y_{j}=x_j-\varepsilon$ for $j=t-2,t-1$ and $y_t=2\varepsilon$, where $0< \varepsilon <c$.

Since $a>b>c$ and the worst case in $G$ is that the number of edges containing the vertex $t$ with weight $2a^2\varepsilon$ is minimum. Then we get that
\begin{eqnarray*}
&&\lambda(G,\vec{y})-\lambda(H_2)\nonumber\\
&=& -2{t-5 \choose 2}a^2\varepsilon-4(t-5)ab\varepsilon-2b^2\varepsilon-2(t-5)ac\varepsilon \nonumber\\
&&+2c^2\varepsilon+8bc\varepsilon+2b^2\varepsilon+4(t-5)ac\varepsilon+4(t-5)ab\varepsilon \nonumber\\
&&+2\left[{t-2 \choose 2}+1-4(t-5)-6\right]a^2\varepsilon+\mathcal{O}(\varepsilon^{2})\nonumber\\
&=&-2(t-3)a^2\varepsilon+2(t-5)ac\varepsilon+8bc\varepsilon+2c^2\varepsilon+\mathcal{O}(\varepsilon^{2})\nonumber\\
&>& \frac{2(t^2-14t+36)a^2\varepsilon}{(t-2)^2}+\mathcal{O}(\varepsilon^{2})\nonumber\\
&>&0,
\end{eqnarray*}
when $\varepsilon$ is small enough. This completes the proof.
\qed

\section{Remarks}
 Results similar to Theorem \ref{theo13} can be obtained for $3$-graphs $G$ containing all but $l$ edges of $[t-1]^{(3)}$ with $l\ge 3$ by modifying the proof of Theorem \ref{theo13}.
   For $r$-graphs, it might be interesting to consider the following conjecture by Peng-Zhao.

\begin{con}{\rm (\cite{PZ})}
Let $t$, $m$ and $r\ge 3$ be positive integers satisfying ${t-1 \choose r} \le m \le {t-1 \choose r} + {t-2 \choose r-1}$.
Let $G$ be an $r$-graph with $m$ edges and $G$ contain a clique of order $t-1$. Then $\lambda(G)=\lambda([t-1]^{(r)})$.
\end{con}

If the above conjecture is true, then combining Theorem \ref{theo3}, one can conclude that for an $r$-graph $G$ with vertex set $[t]$ and $m$ edges containing $[t-1]^{(r)}$, $G$ is dense if and only if $m \ge {t-1 \choose r}+{t-2 \choose r-1}+1$.

Very little is known about dense $r$-graphs. At some point, we thought that an $r$-graph obtained by adding edges to a dense $r$-graph $G$ will be dense. However, this is not true.
\begin{remark}{\rm (Non-monotonicity)}
For a dense $r$-graph $G=(V,E)$, there may exist some $e\in E(G^c)$ such that $H=(V,E\cup \{e\})$ is not dense.
\end{remark}
\noindent{\bf Proof.} For example, let $G$ be a left-compressed $3$-graph on $t$ vertices with $m= {t-1 \choose 3} + {t-2 \choose 2}-1$ edges satisfying $[t-1]^{(3)}\setminus\{(t-3)(t-2)(t-1)\}\subseteq G$ but $[t-1]^{(3)}\nsubseteqq G$ and $(t-3)(t-2)t,(t-3)(t-1)t,(t-2)(t-1)t\notin G$ . Then proposition \ref{prop1} implies that $G$ is dense. Let $H=G\cup\{(t-3)(t-2)(t-1)\}$. It's clear that $[t-1]^{(3)}\subseteq H$ and $|H|= {t-1 \choose 3} + {t-2 \choose 2}$, then $H$ is not dense from Theorem \ref{theorempz}.
\qed

\bigskip
{\bf Acknowledgments.} We thank both reviewers for reading the manuscript carefully, checking all the details  and giving  insightful comments to help improve the manuscript.


\begin{thebibliography}{JluR00}


\bibitem{FF}
P. Frankl, Z. F\"uredi, Extremal problems whose solutions are the
blow-ups of the small Witt-designs, {\it J. Combin. Theory Ser. A.}
{\bf 52}(1989), 129--147.

\bibitem{FR84} P. Frankl and V. R\"{o}dl,  Hypergraphs do not jump, {\it Combinatorica} {\bf 4}(1984), 149-159.

\bibitem{HK}
D. Hefetz and P. Keevash, A hypergraph Tur\'an theorem via lagrangians of intersecting
families, {\it J. Combin. Theory Ser. A} {\bf 120}(2013), 2020--2038.

\bibitem{Keevash} P. Keevash, Hypergrah Tur\'an problems, {\it Surveys in Combinatorics}, Cambridge
University Press, (2011), 83--140.

\bibitem{MS} T.S. Motzkin and E.G. Straus, Maxima for graphs and a new proof of a theorem of Tur\'an, {\it Canad. J. Math} {\bf 17}(1965), 533-540.

\bibitem{PTZ} Y. Peng, Q. Tang and C. Zhao, On Lagrangians of $r$-uniform Hypergraphs, {\it Journal of Combinatorial Optimization} {\bf30(3)}(2015), 812-825.

\bibitem{PZ}  Y. Peng, B. Wu and Y. Yao, A Note on Generalized Lagrangians of Non-uniform Hypergraphs, {\it Order.} (2016), 1--13 (Online).

\bibitem{PZ} Y. Peng, C. Zhao, A Motzkin-Straus type result for $3$-uniform
hypergraphs, {\it Graphs Combin.} {\bf 29}(2013), 681--694.

\bibitem{Sidorenko} A.F. Sidorenko, On the maximal number of edges in a uniform hypergraph that does not contain prohibited subgraphs, {\it Mat. Zametki} {\bf 41}(1987), 247-259.

\bibitem{Sidorenko-89} A.F. Sidorenko, Asymptotic solution for a new class of forbiddenr-graphs, \emph{Combinatorica} {\bf 9}(1989), 207--215.

\bibitem{SPW} Y. Sun, Y. Peng and B. Wu, On Graph-Lagrangians and clique numbers of $3$-uniform hypergraphs, \emph{ Acta Mathematica Sinica, English Series } {\bf 32(8)}(2016), 943--960.

\bibitem{T} J. Talbot, Lagrangians of hypergraphs, {\it Combinatorics, Probability \& Computing} {\bf 11}(2002), 199-216.

\bibitem{TPZZ1}Q. Tang, Y. Peng, X. Zhang and C. Zhao, On graph-Lagrangians of hypergraphs containing dense subgraphs, {\it J. Optim. Theory Appl.} {\bf 163}(2014), no. 1, 31-56.

\bibitem{TPZZ} Q. Tang, Y. Peng, X. Zhang and C. Zhao, Some results on Lagrangians of hypergraphs, {\it Discrete Applied Mathematics} {\bf 166}(2014), 222-238.

\bibitem{TPZZ2} Q. Tang, Y. Peng, X. Zhang and C. Zhao, Connection between the clique number and the Lagrangian of $3$-uniform hypergraphs, {\it Optimization Letters} {\bf 10(4)}(2016), 685-697.



\bibitem{Turan}P. Tur\'an, On an extremal problem in graph theory(in Hungarian), {\it Mat. Fiz. Lapok} {\bf 48}(1941), 436--452.



\end{thebibliography}
\end{document}